\documentclass[12pt]{amsart}
\usepackage{amsmath}
\usepackage{amsfonts,amssymb,amsthm,amstext,amscd,boxedminipage}  %
\usepackage[mathscr]{eucal}

\usepackage{graphicx,psfrag} 
\usepackage{color}
 \numberwithin{equation}{section}

\DeclareMathSymbol{\minus} {\mathord}{operators}{"2D} %

\theoremstyle{plain}
\newtheorem{thm}{Theorem}[section]
\newtheorem{lem}[thm]{Lemma}
\newtheorem{cor}[thm]{Corollary}

\theoremstyle{definition}
\newtheorem{df}[thm]{Definition}
\newtheorem{remark}[thm]{Remark}
\newtheorem{example}[thm]{Example}

\input epsf

\def \H {\mathbb{H}}
\def \Z {\mathbb{Z}}
\def \D {\mathbb{D}}

\def \BRT {{Bollob\'as--Riordan--Tutte }}
\def \IM {{\rm IM}}
\def \sign {{\rm sign}}

\begin{document}

\title[Alternating sum formulae for the determinant]{Alternating sum 
formulae 
for the determinant\\and other link invariants}
\date{\today}

\author[O. Dasbach]{Oliver T. Dasbach}
\address{Department of Mathematics, Louisiana State University,
Baton Rouge, LA 70803}
\email{kasten@math.lsu.edu}

\author[D. Futer]{David Futer}
\address{Department of Mathematics, Michigan State University,
East Lansing, MI 48824}
\email{dfuter@math.msu.edu}

\author[E. Kalfagianni]{Efstratia Kalfagianni}
\address{Department of Mathematics, Michigan State University,
East Lansing, MI 48824}
\email{kalfagia@math.msu.edu}

\author[X.-S. Lin]{Xiao-Song Lin}
\footnote{We regretfully inform you that Xiao-Song Lin  passed away on the 14th of January, 2007.}
 
\author[N. Stoltzfus]{Neal W. Stoltzfus}
\address{Department of Mathematics, Louisiana State University,
Baton Rouge, LA 70803}
\email{stoltz@math.lsu.edu}

\begin{abstract}
A classical result states that the determinant of an alternating link is equal 
to the number of spanning trees in a checkerboard graph of an alternating connected projection 
of the link.

We generalize this result to show that the determinant is the alternating sum of the number of quasi-trees of genus $j$ of the
dessin of a non-alternating link.

Furthermore, we obtain formulas for coefficients of the Jones polynomial by counting quantities on dessins. In particular we will show that the 
$j$-th coefficient of the Jones polynomial is given by sub-dessins of genus less or equal to $j$.
\end{abstract} 

\maketitle

\section{Introduction}

A classical result in knot theory states that the determinant of an alternating link is given by the number of
spanning trees in a checkerboard graph of an alternating, connected link projection (see e.g. \cite{BZ}).
For non-alternating links one has to assign signs to the trees and count the trees with signs, where the geometric meaning of the signs is not apparent.
Ultimately, these theorems are reflected in Kauffman's spanning tree expansion for the Alexander polynomial 
(see \cite{Kauffman:OnKnots, OzsvathSzabo:Alternating})  as well
as Thistlethwaite's spanning tree expansion for the Jones polynomial \cite{Thistlethwaite:SpanningTreeExpansion};  
the determinant is the absolute value of the Alexander polynomial as well as of the Jones polynomial at $-1$.
 
The first purpose of this paper is to show that the determinant theorem for alternating links has a very natural, topological/geometrical generalization to non-alternating links, 
using the framework that we developed in \cite{DFKLS:KauffmanDessins}: Every link diagram induces an embedding of the link into an orientable surface such that the
projection is alternating on that surface. Now the two checkerboard graphs are graphs embedded on surfaces, i.e. dessins d'enfant (aka. combinatorial maps), 
and these two graphs are dual to each other.
The minimal genus of all surfaces coming from that construction is the dessin-genus of the link. 
However, as in \cite{DFKLS:KauffmanDessins} one doesn't need the reference to the surface to construct the dessin directly from the diagram and to compute its genus.
The Jones polynomial can then be considered as an evaluation of the \BRT polynomial \cite{BollobasRiordan:CyclicGraphs} of the dessin \cite{DFKLS:KauffmanDessins}. 
Alternating non-split links are precisely the links of dessin-genus zero. Our determinant formula recovers the classical determinant formula in that case.

For a connected link projection of higher dessin genus we will show that  the determinant
is given as the alternating sum of the number of spanning quasi-trees of genus $j$, as defined below, in the
dessin of the link projection. Thus the sign has a topological/geometrical interpretation in terms of the genus of sub-dessins.  
In particular, we will show that for dessin-genus $1$ projections the determinant is the difference between the number of spanning trees in the
dessin and the number of spanning trees in the dual of the dessin. The class of dessin-genus one knots and links includes
for example all non-alternating pretzel knots.

Every link can be represented as a dessin with one vertex, and we will show that with this representation the numbers of 
$j$-quasi-trees arrise as coefficients of the characteristic polynomial of a certain matrix assigned to the dessin. 
In particular we will obtain a new determinant formula for the determinant of a link which comes solely from the Jones polynomial. 
Recall that the Alexander polynomial - and thus every evaluation of it - can be expressed as a determinant in various ways.
The Jones polynomial, however, is not defined as a determinant. 

The second purpose of the paper is to develop dessin formulas for coefficients of the Jones polynomial.
We will show that  the $j$-th coefficient is completely determined by sub-dessins of genus less or equal to $j$
and we will give formulas for the coefficients.
Again, we will discuss the simplifications in the formulas if the dessin has one vertex.
Starting with the work of the first and fourth author \cite{DL:VolumeIsh}  the coefficients of the Jones polynomial have recently gained a new significance because of their 
relationship to the hyperbolic volume of the link complement. Under certain conditions, the coefficients near the head and the tail of the polynomial give linear upper and lower bounds for the volume. 
In \cite{DL:VolumeIsh, DasbachLin:HeadAndTail} this was done for alternating links and in  \cite{FKP:VolumeJones}  it was generalized to a larger class of links.     

The paper is organized as follows: Section 2 recalls the pertinent results of \cite{DFKLS:KauffmanDessins}. 
In Section 3 we develop the alternating sum formula for the determinant of the link.
Section 4 shows a duality result for
quasi-trees and its application to knots of dessin genus one. 
In Section 5 we look at the situation when the dessin has one vertex. 
Section 6 shows results on the coefficents of the Jones polynomial within the framework of dessins.

{\bf Acknowledgement: } We thank James Oxley for helpful discussions.
The first author was supported in part by NSF grants DMS-0306774 and DMS-0456275 (FRG),  
the second author by NSF grant DMS-0353717 (RTG),
the third author by NSF grants DMS-0306995 and DMS-0456155 (FRG) and
the fifth author by NSF grant DMS-0456275 (FRG).

\section{The Dessin d'enfant coming from a link diagram} \label{sec:dessin}

We recall the basic definitions of \cite{DFKLS:KauffmanDessins}:

A \textit{dessin d'enfant} (combinatorial map, oriented ribbon map) can be viewed
as a multi-graph (i.e. loops and multiple edges are allowed) equipped with a
cyclic order on the edges at every vertex. Isomorphisms between dessins are
graph isomorphisms that preserve the given cyclic order of the edges. 

Equivalently, dessins correspond to graphs embedded on
an orientable surface such that every region in the complement of the graph is a disk.
We call the regions the \emph{faces} of the dessins. 
Thus the genus $g(\D)$ of a dessin $\D$ with $k$ components is determined
by its Euler characteristic: $$\chi(\D)=v(\D)-e(\D) +f(\D) = 2k-2g(\D).$$  

\begin{figure}[htbp]
   \centering
   \includegraphics[width=2.5in]{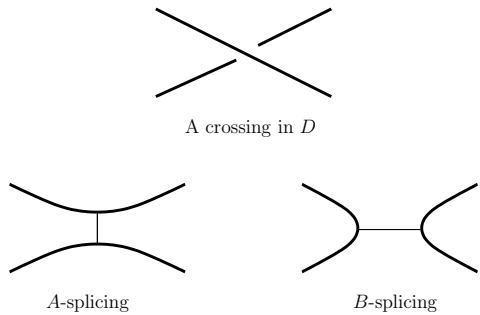} 
   \caption{Splicings of a crossing, $A$-graph and $B$-graph.}
\label{fig:AB-splicing}
\end{figure}

For each Kauffman state of a (connected) link diagram a dessin is constructed as follows: Given a
link diagram $P(K)$ of a link $K$ we have, as in Figure \ref{fig:AB-splicing},
an $A$-splicing and a $B$-splicing at every crossing. For any state assignment
of an $A$ or $B$ at each crossing we obtain a collection of non-intersecting
circles in the plane, together with embedded arcs that record the crossing
splice.  Again, Figure \ref{fig:AB-splicing} shows this situation locally. In
particular, we will consider the state where all splicings are $A$-splicings. 
The collection of circles will be the set of vertices of the dessin.

\medskip

To define the desired dessin associated to a link diagram, we need to define an
orientation on each of the circles resulting from the $A$ or $B$ splicings,
according to a given state assignment. 
We orient the set of circles in the plane by orienting each component clockwise or anti-clockwise according
to whether the circle is inside an odd or even number of circles, respectively.
Given a state assignment $s:E\rightarrow\{A,B\}$ on the crossings (the eventual
edge set $E(\D)$ of the dessin), the associated dessin is constructed by first
resolving all the crossings according to the assigned states and then orienting
the resulting circles according to a given orientation of the plane. 

The vertices of the dessin correspond to the collection of circles and the edges of the
dessin correspond to the crossings. The orientation of the circles defines the
orientation of the edges around the vertices. We will denote the dessin associated to state $s$ by $\D(s)$.
Of particular interest for us will be the dessins $\D(A)$ and $\D(B)$ coming from the states with
all-$A$ splicings and all-$B$ splicings.
For alternating projections of alternating links $\D(A)$ and $\D(B)$ are the two checkerboard graphs
of the link projection. In general, we showed in \cite{DFKLS:KauffmanDessins} that $\D(A)$ and $\D(B)$
are dual to each other.

We will need several different combinatorial measurements of the dessin:
 
\begin{df}\label{def:dessin-counts}
Denote by $v(\D), e(\D)$ and $f(\D)$ the number of vertices, edges and faces of a dessin $\D$.
Furthermore, we define the following quantities:
\begin{eqnarray*}
k(\D) &=& \mbox{the number of connected components of } \D,\\
g(\D) &=& \frac{2k(\D)-v(\D) + e(\D) - f(\D)}{2}, \mbox{the \emph{genus} of } \D,\\
n(\D) &=& e(\D) - v(\D) + k(\D), \mbox{the \emph{nullity} of } \D.
\end{eqnarray*}
\end{df}

The following spanning sub-dessin
expansion was obtained in \cite{DFKLS:KauffmanDessins} by using results of 
\cite{BollobasRiordan:NonOrientableSurfaces}. A \emph{spanning sub-dessin} is obtained from the dessin by
deleting edges. Thus, it has the same vertex set as the dessin.

\begin{thm}\label{thm:sub-dessin}
Let $\langle P \rangle \in \Z[A,A^{-1}]$ be the Kauffman bracket of a connected link projection diagram $P$ and $\D := \D(A)$ be the dessin of $P$ associated to the all-$A$-splicing.
The Kauffman bracket can be computed by the following spanning sub-dessin $\H$ expansion:
\[A^{-e(\D)} \langle P \rangle = A^{2-2v(\D)} (X-1)^{-k(\D)}\sum_{\H \subset \D} (X-1)^{k(\H)}Y^{n(\H)} Z^{g(\H)}
\]
\noindent under the following specialization:  $\{X\rightarrow -A^4,Y\rightarrow A^{-2} \delta,Z\rightarrow  \delta^{-2}\}$
where $\delta:=(-A^2-A^{-2})$.
\end{thm} 
 
\section{Dessins determine the determinant}
 
The determinant of a link is  ubiquitous in knot theory.
It is the absolute value of the Alexander polynomial at -1 as well as the Jones polynomial
at -1. Furthermore, it is the order of the first homology group of the double branched cover
of the link complement. For other interpretations, see e.g. \cite{BZ}.

We find the following definition helpful:

\begin{df} \label{def:quasi-tree}
Let $\D$ be a connected dessin that embeds into the surface $S$.
A spanning quasi-tree of genus $j$ or spanning $j$-quasi-tree in $\D$ is a sub-dessin $\H$ of $\D$ with
$v(\H)$ vertices and $e(\H)$ edges such that $\H$ is connected and spanning and

\begin{enumerate}
\item $\H$ is of genus $j$.
\item $S-\H$ has one component, i.e. $f(\H)=1$.
\item $H$ has $e(\H)=v(\H)-1 + 2 j$ edges.
\end{enumerate}

In particular the spanning $0$-quasi-trees are the regular spanning trees of the graph.
Note that by Definition \ref{def:dessin-counts} either two of the three conditions in Definition \ref{def:quasi-tree}
imply the third one. 
\end{df} 
 
Theorem \ref{thm:sub-dessin} now leads to the following formula for the determinant $\det (K)$ of a link $K$:

\begin{thm} \label{thm:determinant}
Let $P$ be a connected projection of the link $K$ and $\D:=\D(A)$ be the dessin of $P$ associated to the all-$A$ splicing.
Suppose $\D$ is of genus $g(\D)$.

Furthermore, let $s(j,\D)$ be the number of spanning $j$-quasi-trees of $\D$.

Then
$$\det(K)= \left |\sum_{j=0}^{g(\D)} (-1)^j \, s(j,\D) \right |.$$ 
\end{thm} 

\begin{proof}
Recall that the Jones polynomial $J_K(t)$ can be obtained from the Kauffman bracket, 
up to a sign and a power of $t$, by the substitution $t:=A^{-4}$.

By Theorem \ref{thm:sub-dessin} we have for some power $u=u(\D)$:
\begin{align}
\pm J_K(A^{-4})&= A^u \sum_{\H \subset \D} (X-1)^{k(\H)-1} Y^{n(\H)} Z^{g(\H)} \nonumber\\
&= A^u \sum_{\H \subset \D} A^{-2 -2 e(\H)+ 2 v(\H)} \delta^{f(\H)-1}
\end{align}

We are interested in the absolute value of $J_K(-1)$. Thus, $\delta =0$ and, since $k(\H)\leq f(\H)$:
\begin{align}
\left |J_K(-1)\right | &= \left | \sum_{\H \subset \D, f(\H)=1} A^{-2 -2 e(\H) + 2 v(\H)} \right |\\
&= \left | \sum_{\H \subset \D, f(\H)=1} A^{-4 g(\H)} \right |
\end{align}

Collecting the terms of the same genus and setting $A^{-4}:=-1$ proves the claim.
\end{proof}
 
\begin{remark}
For genus $j=0$ we have $s(0,\D)$ is the number of spanning trees in the dessin $\D$. 
Recall that a link has dessin-genus zero if and only it is alternating. Thus, in particular,
we recover the well-known theorem that for alternating links the determinant of a link is
the number of spanning trees in a checkerboard graph of an alternating connected projection.

Theorem \ref{thm:determinant} is a natural generalization of this theorem for non-alternating link
projections. 
\end{remark} 
 
\begin{example} \label{Example821}
 Figure \ref{fig:Eight21p} shows the non-alternating 8-crossing knot $8_{21}$, as drawn by Knotscape (http://www.math.utk.edu/$\sim \!$ morwen/knotscape.html), and 
Figure \ref{fig:Eight21D}   the
all-$A$ associated dessin.

The dessin in Figure \ref{fig:Eight21D} contains 9 spanning trees. Therefore, $s(0,\D)=9$.
A spanning sub-dessin of genus one with $4$ edges must contain either of the two loops and three additional edges.
A simple count yields $24$ of these and thus the determinant of the knot is $24-9=15$.



\begin{figure}[htbp] %
   \centering
   \includegraphics[width=3in, angle=180]{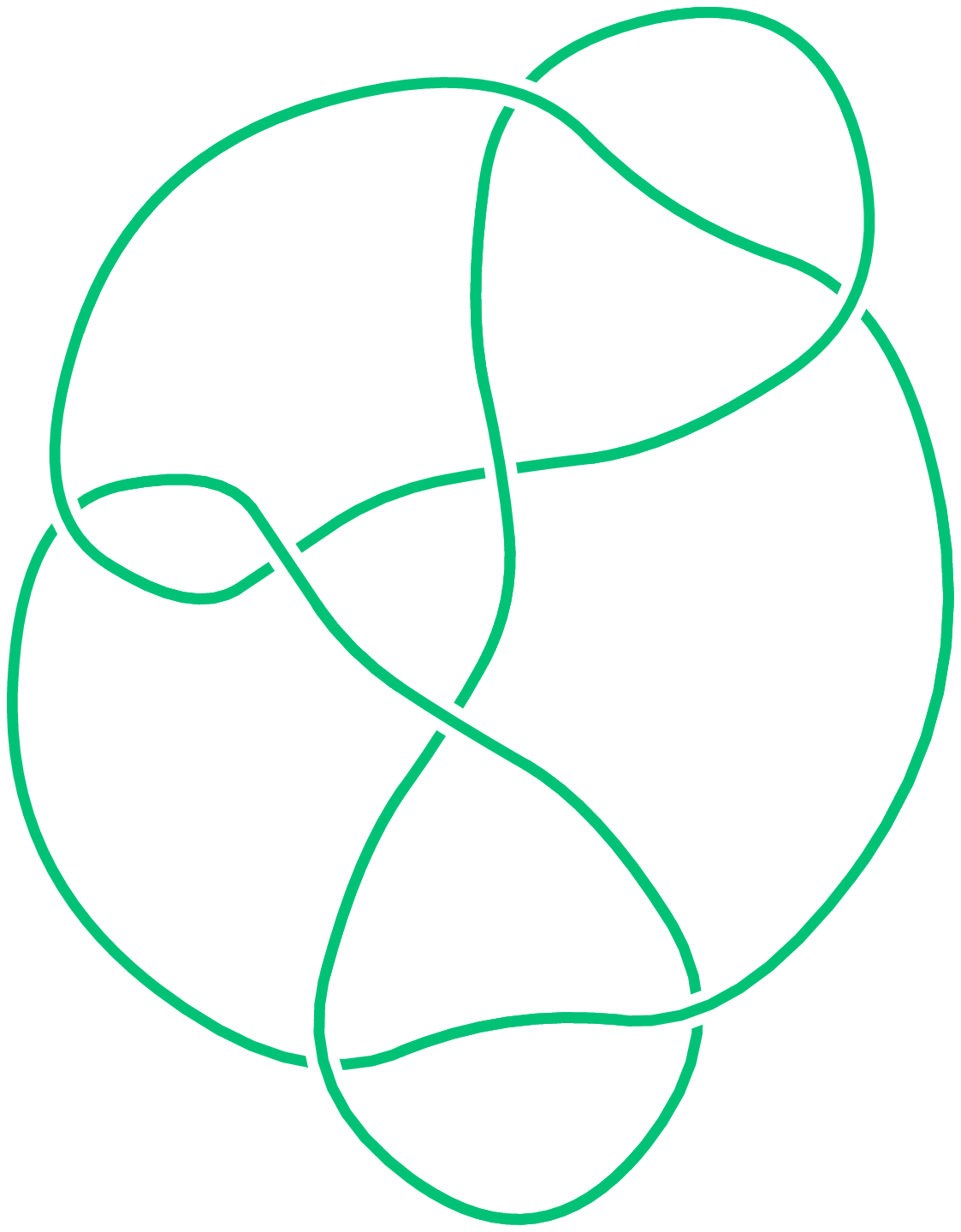} 
   
   \includegraphics[width=2in]{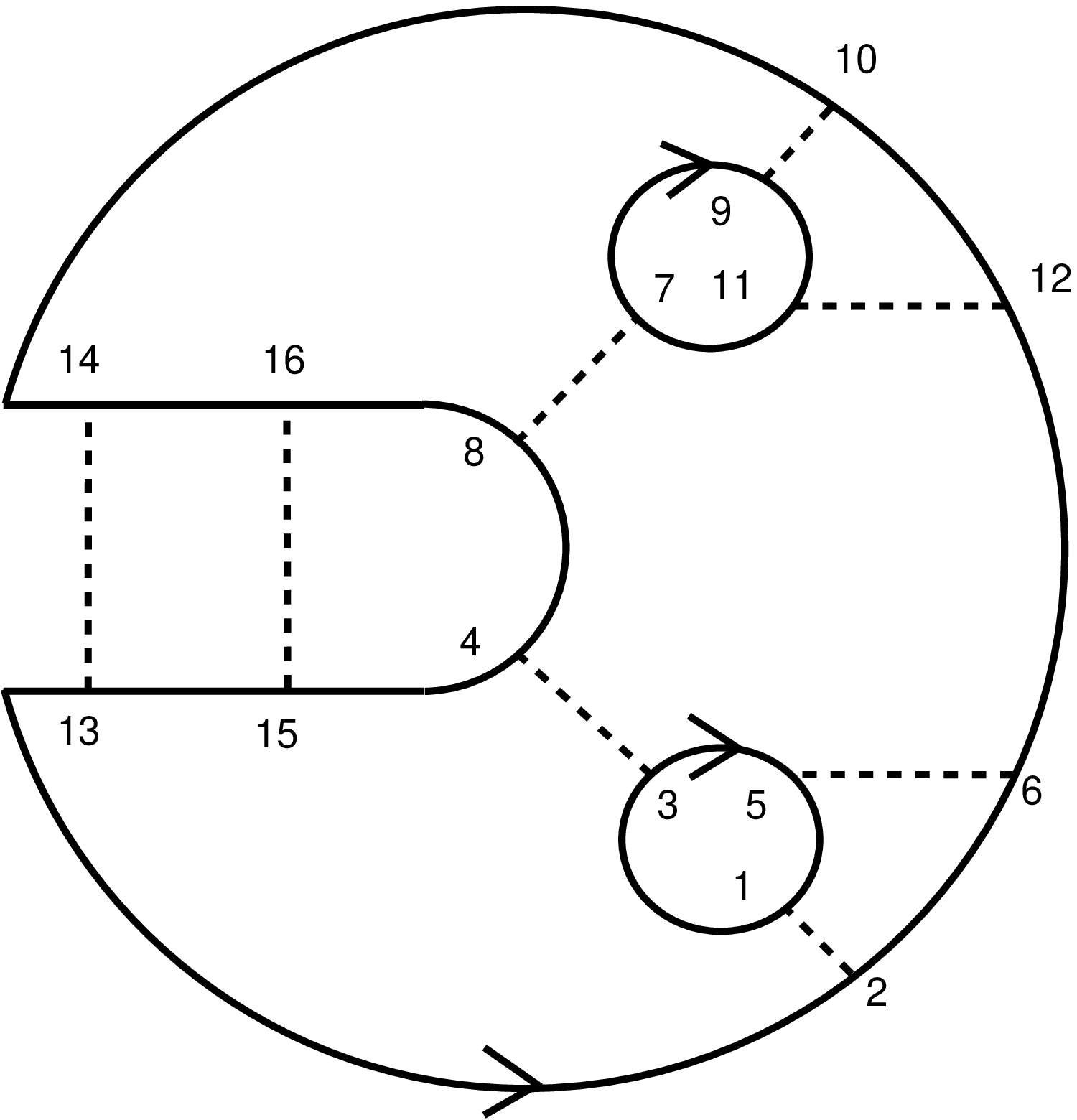} 
   \caption{The eight-crossing knot $8_{21}$ with its all-$A$ splicing projection diagram.}
   	   \label{fig:Eight21p}
\end{figure}

\begin{figure}[htbp] %
   \centering
   \includegraphics[width=2in]{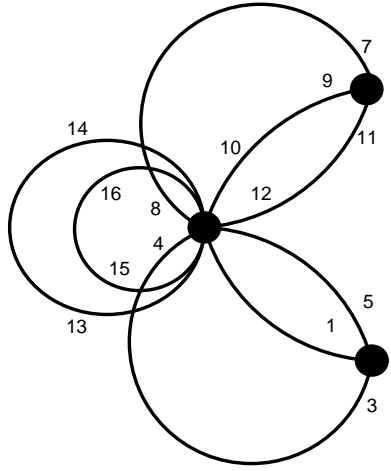} 
   \caption{All-$A$ splicing dessin for $8_{21}$.}
   \label{fig:Eight21D}
\end{figure}
 \end{example}
\medskip

\section{Duality}

The following theorem is a generalization of the result that for planar graphs the spanning trees are in one-one correspondence to the spanning trees of the
dual graphs:

\begin{thm}
Let $\D=\D(A)$ be the dessin of all-$A$ splicings of a connected link projection of a link $L$. Suppose $\D$ is of genus $g(\D)$ and $\D^*$ is the dual of $\D$.

We have: The $j$-quasi-trees of $\D$ are in one-one correspondence to the $(g(\D)-j)$-quasi-trees of $\D^*$. Thus
$$s(j,\D)=s(g(\D)-j,\D^*).$$
\end{thm}

\begin{proof}
Let $\H$ be a spanning $j$-quasi-tree in $\D$. Denote by $\D-\H$ the sub-dessin of $\D$ obtained by removing the edges
of $\H$ from $\D$. From $f(\H)=1$ it follows that the dual $(\D-\H)^*$ is connected and
spanning. Furthermore, $f((\D-\H)^*)=1$.

We have:
\begin{eqnarray*}
v(\H)-e(\H)+f(\H)&=&v(\D)-e(\H)+1 = 2 - 2 j\\
v(\D)-e(\D)+f(\D)&=& 2 - 2 g(\D)
\end{eqnarray*}

Thus,
\begin{eqnarray*}
v((\D-\H)^*)- e((\D-\H)^*)+f((\D-\H)^*) &=& f(\D) - (e(\D)-e(\H)) +1\\
&=& 2 - 2 g(\D) - v(\D) + e(\H) +1\\
&=& 2 - 2(g(\D)-j).
\end{eqnarray*}

Hence, $(\D-\H)^*$ is a $(g(\D)-j)$-quasi-tree in $\D^*$.

\end{proof}

Recall that the dessin-genus zero links are precisely the alternating links. The following corollary generalizes to the class of dessin-genus one links the aforementioned, classical interpretation of the determinant of connected alternating links as the number of spanning trees in its checkerboard graph:
 
\begin{cor} \label{cor:Dessingenus1}
Let $\D=\D(A)$ be the all-$A$ dessin of a connected link projection of a link $L$ and $\D^*$ its dual. Suppose $\D$ is of dessin genus one.
Then

$$
\det (L) =| \# \{\mbox{spanning trees in } \D\} - \# \{\mbox {spanning trees in } \D^*\}|.
$$
\end{cor}

\begin{figure}[htbp] %
   \centering
   \includegraphics[width=2.5in]{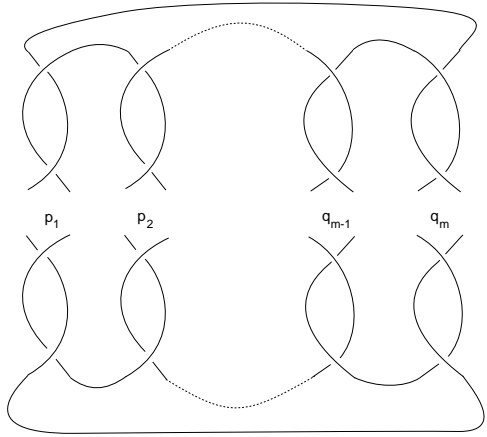} 
   \caption{The $K(p_1, \dots, p_n, -q_1, \dots, -q_m)$ pretzel link.}
   \label{fig:pretzel}
\end{figure}

We apply Corollary \ref{cor:Dessingenus1} to compute the determinants of non-alternating pretzel links.
The Alexander polynomial as well as the Jones polynomial, and consequently the determinant is invariant under mutations
(see e.g. \cite{Lickorish:Book}). Hence, it is sufficient to consider the case of $K(p_1, \dots, p_n, -q_1, \dots, -q_m)$
pretzel links, as depicted in Figure \ref{fig:pretzel}. We assume that the links are non-alternating, i.e. $n\geq 1$ and $m \geq 1$.

\begin{figure}[htbp] %
   \centering
   \includegraphics[width=2.5in]{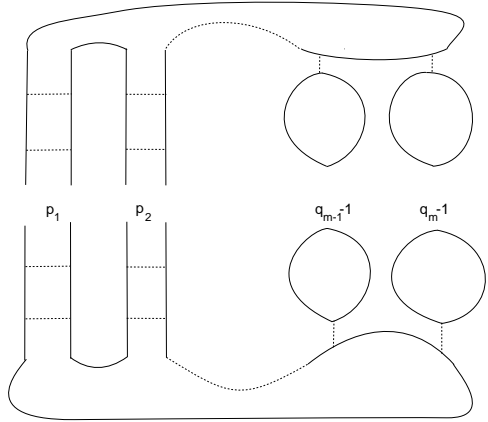} 
   \caption{The all-$A$ splicings of the $K(p_1, \dots, p_n, -q_1, \dots, -q_m)$ pretzel link.}
   \label{fig:asplicedpretzel}
\end{figure}

\begin{lem} Consider the pretzel link $K(p_1,\dots,p_n,-q_1,\dots,-q_m)$, where  $n \geq 1, m \geq 1$ and
$p_i, q_i>0$ for all $i$.  The determinant of $K(p_1,\dots,p_n,-q_1,\dots,-q_m)$,  is
$$
\det(K(p_1, \dots, p_n, -q_1, \dots, -q_m))= \left | \prod_{i=1}^n p_i \prod_{j=1}^m q_j \left(\sum_{i=1}^n \frac 1 {p_i} - \sum_{j=1}^m \frac 1 {q_j} \right ) \right |
$$
\end{lem}

\begin{proof}
Figure \ref{fig:asplicedpretzel} shows the all-$A$ splicing diagram of these links.
The all-$A$ dessin $\D = \D(A)$ has 
$$v(\D)=n+\sum_{j=1}^m (q_j-1)= n-m + \sum_{j=1}^m q_j$$
vertices and $e(\D)=\sum_{i=1}^n p_i + \sum_{j=1}^m q_j$ edges.
For the numbers of faces we have to count the vertices in the all-$B$ dessin.
We compute:
$$f(\D)= m-n+\sum_{i=1}^n p_i.$$
Now we get for the Euler characteristic:
$$\chi(\D)=v(\D)-e(\D)+f(\D)=0$$
and thus the dessin-genus is one.

It remains to compute the difference between the number of spanning trees in the dessin and the number of spanning trees in its dual.
This is a simple counting argument.
\end{proof}

\begin{remark}
The class of dessin-genus one knots and links is quite rich. 
For  example, it contains all non-alternating Montesinos links. It also  
contains all semi-alternating links (whose diagrams are constructed  
by joining together two alternating tangles, and thus have exactly  
two over-over crossing arcs and two under-under arcs).
\end{remark}

\section{Dessins with one vertex}

\subsection{Link projection modifications} 

Here we show that every knot/link admits a projection with respect to which the all-$A$ dessin has one vertex.
Such dessins are useful for computations.

\begin{lem}\label{Onevertexdessin} Let $\tilde P$ be a projection
of a link $L$ with corresponding  all-$A$ dessin $\tilde \D$.  Then $\tilde P$ can be modified by
Reidemeister moves  to a new  a projection $P$ such that the corresponding dessin $\D=\D(A)$ has one vertex.
Furthermore, we have:
\begin{enumerate}
\item $e(\tilde \D)+ 2 v(\tilde D)-2 = e (\D)$
\item $g(\tilde \D)+ v(\tilde D) -1 = g(\D).$
\end{enumerate}
\end{lem}

\begin{proof}
For a connected projection of the link $L$ consider the collection of circles that we obtain by an all-$A$ splicing of the crossings.
If there is only one circle we are done. Otherwise, one can perform a Reidemeister move II near a crossing on two arcs that lie on two neighbor circles as in Figure \ref{fig:vertexredux}.

\begin{figure}[htbp] %
   \centering
   \includegraphics[width=3.5in]{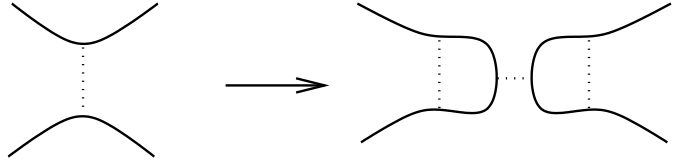} 
   \caption{Reduction of the number of vertices by a Reidemeister II move}
   \label{fig:vertexredux}
\end{figure}

The new projection will have one circle less in its all-$A$ splicing diagram. Also two crossings were added and a new face was created. 
If the link projection is  non-connected one can transform it by Reidemeister II moves into a connected link projection.  It is easy to check that the genus behaves as predicted.
The claim follows.
\end{proof} 

\begin{remark}
Dessins with one vertex are equivalent to Manturov's ``d-diagrams" \cite{Manturov:Ddiagrams}. Note that the procedure of using just Reidemeister moves of type II
is similar in spirit to Vogel's proof of the Alexander theorem \cite{Vogel:AlexanderTheorem, BirmanBrendle:BraidSurvey}.
\end{remark}

\subsection{The determinant of dessins with one vertex}

Dessins with one vertex can also described as chord diagrams.
The circle of the chord diagram corresponds to the vertex and the chords correspond to the edges. 
In our construction the circle of the chord diagram is the unique circle of the state resolution, and the chords correspond to the crossings.
The cyclic orientation at the vertex induces the order of the chords around the circle.
For each chord diagram $\D$ one can assign an intersection matrix \cite{CDL, BN-G:Melvin-Morton} as follows:
Fix a base point on the circle, disjoint from the chords and number the chords consecutively.

The intersection matrix is given by:

$$\IM(\D)_{i j}= \left \{ \begin{array}{ll}
\sign (i-j) & \mbox{if the $i$-th chord and the $j$-th chord intersect} \\
0 & \mbox {else}
\end{array} \right .$$

Recall that the number of spanning $j$-quasi trees in $\D$ was denoted by $s(j,D)$.
Now:

\begin{thm} 
For a dessin $\D$ with one vertex the characteristic polynomial of $\, \IM(\D)$ satisfies:

$$\det(\IM(\D) - x I) = (-1)^{m} \sum_{j=0}^{\lfloor \frac{m} 2 \rfloor } s(j,\D) x^{m-2 j},$$  
where $m=e(\D)$ is the number of chords, i.e. the number of edges in the dessin.

In particular

$$\det(\D) = | \det (\IM(\D) - \sqrt {-1} I)|.$$

\end{thm}

\begin{proof}
The result follows from combining Theorem \ref{thm:determinant} and a result of Bar-Natan and Garoufalidis \cite{BN-G:Melvin-Morton}.
Bar-Natan and Garoufalidis use chord diagrams to study weight systems coming from Vassiliev invariant theory, thus in a different setting than we do.
However, by \cite{BN-G:Melvin-Morton} for a chord diagram $\D$ the determinant $\det(\IM(\D))$ is either $0$ or $1$ and, translated in our language, it is $1$ precisely if $f(\D)=1.$

Furthermore, since $f(\D)-1$ and the number of edges have the same parity, we know that $\det(\IM(\D))=0$ for an odd number of edges.
 
The matrix $\IM(\D)$ has zeroes on the diagonal.
Thus the coefficient of $x^{m-j}$ in $\det(\IM(\D) - x I)$ is $(-1)^{m-j}$ the sum over the determinants of all $j \times j$ submatrices that are obtained by deleting $m-j$ rows and the $m-j$ corresponding columns in the matrix $\IM(\D)$. Those submatrices are precisely $\IM(\H)$ for $\H$ a subdessin of $\D$ with $j$ edges. In particular the determinant of $\IM(\H)$ is zero for $j$ odd. For $j$ even we know that $\det(\IM(\H)))=1$ if $f(\H)=1$ and $0$ otherwise. Since for $1$-vertex dessins $D$ the genus $2 g(\D)= e(\D)-f(\D)+1$ those $\H$ with $f(\H)$ are precisely the $j$-quasi-trees.
This, together with Theorem \ref{thm:determinant}  implies the claim.
\end{proof}

\begin{example}
The $(p,q)$-twist knots as in Figure \ref{fig:pq-twist} have an all-$A$-dessin with one vertex.

\begin{figure}[htbp] %
   \centering
   \includegraphics[width=4in]{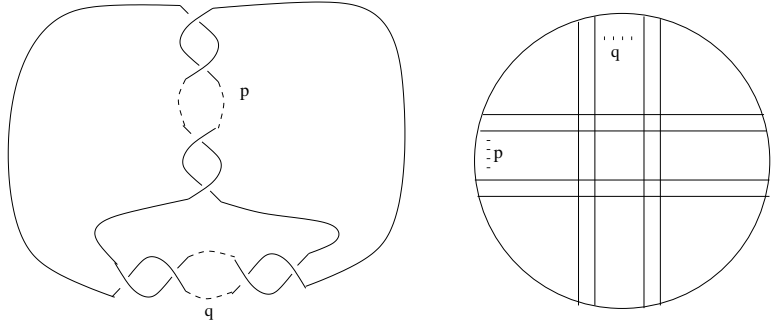} 
   \caption{The $(p,q)$-twist knot and its all-$A$ splicing dessin in chord diagram form.}
   \label{fig:pq-twist}
\end{figure}

The figure-8 knot is given as the $(2,3)$-twist knot. Its intersection matrix is
$$\IM(\D) = \left ( \begin{array}{ccccc}
0&0&-1&-1&-1\\
0&0&-1&-1&-1\\
1&1&0&0&0\\
1&1&0&0&0\\
1&1&0&0&0.
\end{array} \right )
$$
The characteristic polynomial of $\IM(\D)$ is $-6 x^3 - x^5.$ In particular, the determinant of the figure-8 knot is $6-1=5$.
\end{example}

\subsection{The Jones polynomial at $t=-2$.}

By work of Jaeger, Vertigan and Welsh \cite{JVW} evaluating the Jones polynomial is $\#P$-hard at all points, except at
eight points: All fourth and sixth roots of unity. In particular, the determinant arises as one of these exceptional points.
However, letting computational complexity aside, Theorem \ref{thm:sub-dessin} gives an interesting
formula in terms of the genus for yet another point: $t=A^{-4}=-2$:

\begin{lem}\label{Jonesat-2}
Let $P$ be the projection of a link $K$ with dessin $\D=\D(A)$ such that $\D$ has one vertex. Then the Kauffman bracket at
$t=A^{-4}:=-2$ evaluates to

$$\langle P \rangle = A^{e(\D)} \sum_{\H \subset \D} (A^{-4})^{g(\H)}.
$$.
\end{lem}

\begin{proof}
Theorem \ref{thm:sub-dessin}, after substitution, yields the following sub-dessin expansion for the Kauffman bracket of $P$:
$$\langle P \rangle = \sum_{\H \subset \D} A^{e(\D)- 2 e(\H)} \left ( -A^2-A^{-2} \right )^ {f(\H)-1}.$$
The term
$$A^{-2} (-A^2-A^{-2})=(-1-A^{-4})$$
is $1$ at $t=A^{-4}=-2$ and, with $v(\D)=v(\H)$ for all spanning sub-dessin $\H$ of $\D$, the claim follows.
\end{proof}
 
\section{Dessins and the coefficients of the Jones Polynomial}

Let $P$ be a connected  projection of a link $L$, with corresponding 
all-$A$ dessin 
 $\D:=\D(A)$ and let
 
\begin{equation}\langle P \rangle = \sum_{\H \subset \D} A^{e(\D)-2e(\H)}   (-A^2 -A^{-2})^{f(\H)-1}   \label{Eq:KauffmanBracket}  \end{equation}

\noindent
denote the spanning sub-dessin expansion of the Kauffman bracket of $P$ as obtained earlier.
Let ${\H}_0\subset {\D}$ denote the spanning sub-dessin that contains no edges (so
$f({\H}_0)=v(\D)$ and $e({\H}_0)=0$) and let $M:=M(P)$ and $m:=m(P)$ denote
the maximum and minimum powers of $A$ that occur in the terms that lead to $\langle P \rangle$.
We have 
$$M(P)\leq e(\D)+2v(\D)-2,$$ 
and the exponent $e(\D)+2v(\D)-2$
is realized by ${\H}_0$; see Lemma 7.1, \cite{DFKLS:KauffmanDessins}. Let
$a_M$ denote the coefficient  of the extreme term $A^{e(\D)+2v(\D)-2} $ of $ \langle P \rangle$. Below we will give formulae
for $a_M$; 
similar formulae can be obtained for the lowest coefficient, say $a_m$, if one replaces the 
the all-$A$ dessin with the all-$B$ dessin in the statements below.
We should note that $a_M$ is not, in general, the first non-vanishing coefficient of the Jones polynomial of
$L$. Indeed, the exponent ${e(\D)+2v(\D)-2}$ as well as the expression for $a_M$ we obtain below,
depends on the projection $P$ and it is not, in general, an invariant of $L$. In particular, $a_M$ might be zero and, for example, we will show that this is the case 
in Example \ref{example:821}.

The following theorem extends and recovers results of Bae and Morton, and Manch{\'o}n \cite{MortonBae:ExtremeTerms, Manchon:ExtremeCoefficients} within the dessin framework.
\begin{thm} \label{formula} We have 
\begin{enumerate} 
\item For $l\geq 0$ let $a_{M-l}$ denote the coefficient of $A^{e(\D)+2 v(\D)-2 -4l}$ in the Kauffman bracket $\langle P \rangle$.
Then, the term $a_{M-l}$ only depends on spanning sub-dessins $\H \subset \D$ of genus  $g(\H)\leq l$.

\item The highest term is given by
\begin{equation} a_M = \sum_{\H \subset \D, \, g(\H)=0=k(\H)-v(\D)} (-1)^{v(\D)+e(\H)-1}. \end{equation}
\noindent
In particular, if $\D$ does not contain any loops then $a_M=(-1)^{v(\D)-1}$ and the only contribution 
comes  from ${\H}_0$.

\end{enumerate}
\end{thm}

\begin{proof}
The contribution of a spanning $\H \subset \D$ to $ \langle P \rangle$ is
given by
 
\begin{equation}
X_{\H}:= 
A^{e(\D)-2e(\H)}   (-A^2 -A^{-2})^{f(\H)-1}. 
\end{equation}
\noindent
A typical monomial of $X_{\H}$ is of the form
$A^{e(\D)-2e(\H)+2f(\H)-2-4s}$, for $0\leq s \leq f(\H)-1$.
For a monomial to contribute to $a_{M-l}$ we must have
\begin{equation}
e(\D)-2e(\H)+2f(\H)-2-4s=e(\D)+2v(\D)-2-4l,
\end{equation}
or 
\begin{equation}
 f(\H)=v(\D)+e(\H)+2s-2l,
\end{equation}

Now we have
\begin{eqnarray*}
2g(\H)&=&2k(\H)-v(\D) + e(\H) - f(\H)\\
&=& 2k(\H)-2v(\D) + 2l-2s,\\
\end{eqnarray*}
or  $g(\H)=k(\H)-v(\D)+ l-s$.
But since $v(\D)\geq k(\H)$
(every component must have a vertex) and $s \geq 0$
we conclude that 

$$l = g(\H)+ v(\D)-k(\H) + s \geq g(\H),$$
as desired.
Now to get the claims for $a_M$: Note that
for a monomial of $X_{\H}$ to contribute to $a_M$ we must have
\begin{equation}
 g(\H)=k(\H)-v(\D) -s
\end{equation}
which implies
that $s=g(\H)=0$ and $v(\D)= k(\H)$.
It follows that $\H$ contributes to $a_M$ if and only if all of the following conditions are satisfied:

\begin{enumerate}
\item $f(\H)=v(\D)+e(\H)$.
\item $k(\H)=v(\D)$. Thus $\H$ consists of $k:=k(\H)$ components each of which has exactly one vertex
and  either $\H$ has no edges or every edge is a loop.
\item $g(\H)=0$.
\item  the contribution of $\H$ to $a_M$ is $(-1)^{f(\H)-1}$. 
\end{enumerate}

This finishes the proof of the theorem. \end{proof}

\begin{example} \label{example:821}  The all-$A$ dessin of Figure \ref{fig:Eight21D} contains one subdessin
with no edges, two subdessins with exactly one loop and one subdessin of genus zero with two loops.
Thus $a_M=0$. 
\end{example}

A connected link projection is called $A$-adequate iff the all-$A$ dessin $\D(A)$ contains no loops; alternating links admit such projections.  We consider two edges as equivalent if they connect the same two vertices. Let  $e'=e'(\D(A)$ denote the number of edges of equivalence classes of edges.

The following is an extension in \cite{Stoimenow:SecondCoefficient}  to the class of adequate links of a result in \cite{DL:VolumeIsh}  for alternating links.
We will give the dessin proof for completeness, since it shows a subtlety when dealing with dessins in our context: 
Not all dessins can occur as a dessin of a link diagram.

\begin{cor} 
For $A$-adequate diagrams $a_{M-1}= (-1)^{v} (e'-v+1)$
 \end{cor}

\begin{proof} With the notation and setting of the proof of Theorem
\ref{formula} we are looking to calculate the coefficient of the power
$A^{e(\D)+2 v(\D)-6}$. The analysis in the proof of Theorem \ref{formula}
implies that a spanning sub-dessin $\H \subset \D$ contributes to
$a_{M-1}$ if it satisfies one of the following:

(1) $v(\H)=k(\H)$ and $g(\H)=1$.

(2) $v(\H)=k(\H)$ and $g(\H)=0$.

(3) $v(\H)=k(\H)+1$ and $g(\H)=0$.

Since the link is adequate  $\D(A)$  contains no loops and we cannot have any $\H$ as in (1).
Furthermore, the only $\H$ with the properties of (2) is the sub-dessin $\H_0$ 
that contains no edges. Finally the only case that occurs in (3) consists of those sub-dessins
$\H_1$ that are obtained from $\H _0$ by adding edges
between a  pair of vertices. The dessin is special since it comes from a link diagram. Each vertex in the dessin represents a circle in the all-$A$ splicing diagram of the link and
each edge represents an edge there.  Because these edges do not intersect  $\H_1$ must have genus $0$.

Note that any sub-dessin  $\H'\subset \H_1$ is either $\H_0$ or is of the sort described in (3). 
We will call  $\H_1$ maximal if its not properly contained in one of the same type with more edges.
Thus there are $e'$ maximal $\H_1$ for $\D(A)$. The contribution of $\H_1$ to $a_{M-1}$ is $(-1)^{v(\D)-3+e(\H_1)}$.
Thus the contribution of all $\H' \subset \H_1$ that are not $\H_0$ is
$$\sum_{j=1}^{e(\H_1)} {e(\H_1) \choose j}   (-1)^{v(\D)-3+j} = (-1)^{v(\D)}.$$

Thus, the total contribution in $a_{M-1}$ of all such terms is $(-1)^v e'$. 

To finish the proof, observe that the contribution of $\H_0$ comes from the second term of
the binomial expansion 
\begin{equation}
X_{\H_0}:= 
A^{e(\D)}   (-A^2 -A^{-2})^{f(\H_0)-1}. 
\end{equation}
\noindent
Since $f(\H_0)=v$ this later contribution is $(-1)^{v-1} (v-1)$. \end{proof}

The expression in Theorem \ref{formula} becomes simpler, and the lower order terms easier to express, if the dessin $\D$ has only one vertex. By Lemma \ref{Onevertexdessin}
the projection $P$ can always be chosen so that this is the case.

\begin{cor} \label{Cor:CoeffOneVertex}          Suppose $P$ is a connected link projection such that $\D=\D(A)$ has one vertex. Then,
\begin{equation}
a_{M-l} = \sum_{\H \subset \D,\, g(\H)=0}^{g(\H)=l} (-1)^{e(\H)} {e(\H)-2 g(\H) \choose l - g(\H) }.
\end{equation}

In particular,
\begin{equation} a_M = \sum_{\H \subset \D ,\, g(\H)=0} (-1)^{e(\H)} \end{equation}
and
\begin{equation} a_{M-1} = \sum_{\H \subset \D,\, g(\H)=1 } (-1)^{e(\H)} +
                           \sum_{\H \subset \D,\, g(\H)=0 } (-1)^{e(\H)} e(\H)
                           \end{equation}
                  
\end{cor}

\begin{proof}
For a 1-vertex dessin $\D$ we have
$$k(\D)=v(\D)=1 \mbox{ and, thus } 2 g(\D)=e(\D)-f(\D)+1.$$

Now Equation (\ref{Eq:KauffmanBracket}) simplifies to
\begin{eqnarray*}
\langle P \rangle &=& \sum_{\H \subset \D} A^{e(\D)-2 e(\H)+ 2 f(\H)-2} \left (-1 - A^{-4} \right ) ^{f(\H)-1}\\
&=& \sum_{\H \subset \D} A^{e(\D)-4 g(\H)} \left ( -1 - A^{-4} \right )^{e(\H)- 2 g(\H) }.
\end{eqnarray*}
The claim follows from collecting the terms.
\end{proof}

Parallel edges, i.e. neighboring edges that are parallel in the chord diagram,  in a dessin are special since they correspond to twists in the diagram. It is useful to introduce weighted dessins:
Collect all edges, say $\mu(e)-1$ edges parallel to a given edge $c$ and replace this set by $c$ weighted with weight $\mu(c)$.
Note, that $\tilde \D$ has the same genus as $\D$.

\begin{cor} \label{weighted_dessin} For a knot projection with a $1$-vertex dessin $\D$ and weighted dessin $\tilde D$ we have:
$$\langle P \rangle = \sum_{\tilde \H \subset \tilde \D} A^{e(\D)-4 g(\tilde \H)} (-1-A^{-4})^{-2 g(\tilde \H)} 
\prod_{c \in \tilde \H} (-1 - A^{-4 \mu(c)}).$$
\end{cor}

\begin{proof}
For a given edge $c$
collect in
$$\langle P \rangle =\sum_{\H \subset \D} A^{e(\D)-4 g(\H)} \left ( -1 - A^{-4} \right )^{e(\H)- 2 g(\H) }$$
all terms where $\H$ contains an edge parallel to $c$.
This sub-sum is
\begin{eqnarray*}
\sum_{\H \subset \D, \, \H \, \mbox{\tiny contains edge parallel to } c  } A^{e(\D)-4 g(\H)} \left ( -1 - A^{-4} \right )^{e(\H)- 2 g(\H)}&=& \\
\sum_{\H \subset \D, \, \tilde \H=\H -\{ \mbox{\tiny edges parallel to } c \}  \cup {c}}  \sum_{j=1}^{\mu(c)} {\mu(c) \choose j} A^{e(\D)-4 g(\H)} \left ( -1 - A^{-4} \right )^{e(\tilde \H)-1+j- 2 g(\H)}&=&\\
\sum_{\H \subset \D, \, \tilde \H=\H -\{ \mbox{\tiny edges parallel to } c \}  \cup {c}}  A^{e(\D)-4 g(\H)} \left ( -1 - A^{-4} \right )^{e(\tilde \H)-1- 2 g(\H)}  \left ( -1 - A^{-4 \mu(c)} \right ) &&
\end{eqnarray*}
The claim follows by repeating this procedure for each edge $c$.
\end{proof}

\begin{example}
The $(p,q)$-twist knot is represented by the weighted, $1$-vertex dessin with two intersecting edges, one with weight $p$ and one with weight $q$. By Corollary \ref{weighted_dessin} its Kauffman bracket is:
$$
A^{-p-q} \langle P \rangle = 1 + (-1 - A^{-4 p}) + (-1-A^{-4 q}) + A^{-4} (-1 - A^{-4})^{-2} (-1 -A^{-4 p}) (-1 - A^{-4 q}).
$$  
\end{example}

\begin{remark}
Corollary \ref{Cor:CoeffOneVertex} implies the following for the first coefficient $a_M$.
Suppose $\D$ is a $1$-vertex, genus $0$ dessin with at least one edge. Then every sub-dessin also has genus $0$.

Thus,

$$\sum_{\H \subset \D,\, g(\H)=0} (-1)^{e(\H)} = \sum_{j=0}^{e(\D)} {e(\D) \choose j} (-1)^j = 0.$$

For an arbitrary dessin let $\H_1, \dots, \H_n$ be the maximal genus 0 subdessins of $\D$.
Define a function $\phi$ on dessins which is $1$ if the dessin contains no edges and $0$ otherwise.

Then
$$a_M= \sum_{i} \phi (\H_i) - \sum_{i,j,\, i < j} \phi(\H_i \cap \H_j) + \sum_{i,j,k,\, i<j<k} \phi(\H_i \cap \H_j \cap \H_k) - \dots$$
\end{remark}
 
\bibliography{../linklit}
\bibliographystyle {alpha}
 
\end{document}